\documentclass[amstex,12pt]{article}
\usepackage{amsfonts}
\usepackage{amssymb}
\usepackage{amsmath}
\usepackage{mathrsfs}
 \usepackage{indentfirst}

\usepackage{amssymb,amsmath,amsthm}
\newcommand{\beq}{\begin{equation}}
\newcommand{\eeq}{\end{equation}}
\newcommand{\beqn}{\begin{eqnarray}}
\newcommand{\eeqn}{\end{eqnarray}}

\def\e{\hat{\mathbb{E}}}
\def\h{\mathcal {H}}

\def\k{\Omega,\mathcal {H},\hat{\mathbb{E}}}
\renewcommand\o{{\Omega}}

\begin{document}
\small

\title{A weighted central limit theorem under sublinear expectations }
\author{{Defei Zhang\thanks{\small{Corresponding author. E-mail: zhdefei@163.com.}}} \\ \ \it \small{ School of Mathematics,
Shandong University, Jinan 250100, China }}
\date{}
\maketitle
\begin{abstract} In this paper, we investigate a  central limit theorem for weighted sums of independent random variables under sublinear expectations. It is turned out that our results are  natural extensions of the  results obtained by Peng  and Li and Shi.
\end{abstract}
\section{Introduction}
Many useful linear statistics based on a random sample are weighted sums of i.i.d. random variables, for example,  least-squares estimators, nonparametric regression function estimators and jackknife estimates, etc. In this respect, studies of limit theorem for these weighted sums have demonstrated significant progress in probability theory with applications in mathematical statistics, see \cite{c} et al.

In 2006, motivated by measures of risk, super-hedge pricing and model uncertainty in finance, Peng (see \cite{p1}, \cite{p2}) introduced the  notion of independent and identically distributed random variables under sublinear expectations. He also
obtained the law of large numbers and a new central limit  theorem
in the framework of sublinear expectation space in \cite{p3}. The notions of G-distributed, G-Brownian motion and G-expectation were introduced by Peng in \cite{p4}-\cite{p6}. These notions have very rich and interesting new structures which nontrivially generalize the classical ones. Since then, many results have been obtained in [2-9], [11] and [18-21], etc.

However, to our knowledge, there are only few research results about central limit theorem under sublinear expectations. Up to now, Hu and Zhang \cite{h} obtained the central limit theorem for capacities induced by sublinear expectations.  A central limit theorem without the requirement of identical distribution was obtained by  Li and Shi in \cite{l}.  In this paper, we obtain  a  central limit theorem for weighted sums of independent random variables under sublinear expectations. It is turned out that our results are  natural extensions of the results obtained by Peng \cite{p3} and Li and Shi \cite{l}.

This paper is organized as follows: in Section 2, we recall some notions and propositions under sublinear
expectations that are useful in this paper. In Section 3, we give our main result including the proof.
\section{Preliminaries}

In this section, we  recall some basic notions such as sublinear
expectation, identical distribution, independence and G-distributed
(see Peng \cite{p1}-\cite{p6}), et al.

For a given positive integer $n$, we  denote by $\langle
x,y\rangle$ the scalar product of $x, y \in R^n$ and by $|x|
=\langle x,x\rangle^{1/2}$ the Euclidean norm of $x$. We denote by
$S(n)$ the collection of $n\times n$ symmetric matrices.

Let $\Omega$ be a given set and let $\mathcal {H}$ be a linear space
of real-valued  functions defined on $\Omega$ such that if
$X_1,\cdots ,X_n\in \h$, then $\varphi(X_1,\cdots ,X_n)\in \h$ for
each $\varphi \in C_{l,lip}(R^n)$, where $C_{l,lip}(R^n)$ denotes the linear
space of (local Lipschitz) functions $\varphi$ satisfying
$$|\varphi(x)-\varphi(y)|\leq C(1+|x|^m+|y|^m)|x-y|,\forall x,y\in
R^n,$$ for some constants $C > 0; m\in N$ depending on $\varphi.$
$\h$ is considered as a space of  random variables. In this case $X=(X_1,\cdots ,X_n )$ is called an $n$-dimensional vector,
denoted by $X\in \h^n$. we  also denote
\\$\bullet$ $C_{b,lip}(R^n):$ the space of bounded
and Lipschitz continuous functions.\\ \textbf{Definition 2.1.}  A nonlinear expectation $\e$ on ${\cal H}$ is a functional
$\e$ : ${\cal H}\mapsto\mathbb R$ satisfying the following
properties: for all $X$, $Y\in{\cal H}$, we have\\
(a) Monotonicity:\ \ \ \ If $X\geq Y$ then
$\e[X]\geq\e[Y]$.\\
(b) Constant preserving: $\e[c]=c$.\\
The triple $(\Omega,{\cal H},\e)$ is called a nonlinear
expectation space (compare with a probability space $(\Omega,{\cal
F},P)$). We are mainly concerned with sublinear expectation where
the expectation $\e$ satisfies
also\\
(c) Sub-additivity:\ \ \ \ $\e[X]-\e[Y]\leq
\e[X-Y]$.\\
(d)Positive homogeneity: $\e[\lambda
X]=\lambda\e[X]$, $\forall\lambda\geq0$.\\
If only (c) and (d) are satisfied, $\e$ is called a
sublinear functional.

Let $X=(X_1,\cdots ,X_n)$ be a given $n$-dimensional random vector
on a sublinear expectation space $(\o,\h,\e)$. We define a
functional on $C_{l,lip}(R^n)$ by
$$\hat{\mathbb{F}}_X[\varphi]:=\e[\varphi(X)]:\varphi\in
C_{l,lip}(R^n)\rightarrow R. $$ $\hat{\mathbb{F}}_X$ is called the
distribution of $X$. \\
\textbf{Definition 2.2.} Let $X_1$ and $X_2$ be two $n$-dimensional
random vectors defined respectively in sublinear expectation spaces
$(\o_1,\h_1,\e_1)$ and $(\o_2,\h_2,\e_2)$. They are called
identically distributed, denoted by $X_1 \stackrel{d}{=}X_2$,
if$$\e_1[\varphi(X_1)]=\e_2[\varphi(X_2)],\forall \varphi\in
C_{l,lip}(R^n).$$

If the distribution $\hat{\mathbb{F}}_X$ of $X\in\h$ is not a linear
expectation, then $X$ is said to have distributional uncertainty.
The distribution of $X$ has the following four typical parameters:
$$\bar{\mu}:=\e[X];\  \underline{\mu}:=-\e[-X];\  \bar{\sigma}^2:=\e[X^2];\
\underline{\sigma}^2:=-\e[-X^2].$$ The intervals $[\underline{\mu},
\bar{\mu}]$ and $[\underline{\sigma}^2, \bar{\sigma}^2]$
characterize the mean-uncertainty and the variance-uncertainty of
$X$ respectively. \\
\textbf{Definition 2.3.} In a sublinear expectation space $(\k)$, a
random vector $Y\in \h^n$ is said to be independent from another
random vector  $X\in \h^m$  under $\e[\cdot]$ if for each test
function $\varphi\in C_{l,lip}(R^{m+n})$, we have
$$\e[\varphi(X,Y)]=\e[\e[\varphi(x,Y)]_{x=X}].$$

It is important to note that under sublinear expectation the
condition ``$Y$ is independent from $X$'' does not imply automatically
that  ``$X$ is independent from $Y$''.
\\\textbf{Definition 2.4 (maximal distribution).} A $d$-dimensional
random vector $\eta=(\eta_1,\cdots,\eta_d)$ on a sublinear
expectation space $(\k)$ is called maximal distribution if there
exists a bounded, closed and convex subset $\Gamma\subset R^d$ such
that $$\e[\varphi(\eta)]=\max\limits_{y\in\Gamma}\varphi(y),\varphi\in
C_{l,lip}(R^d).$$

 It is
easy to check that this maximal distributed random vector $\eta$
satisfies $$a\eta+b\bar{\eta}\stackrel{d}{=}(a+b)\eta,\ \ \forall a,b\geq 0,$$
where $\bar{\eta}$ is an independent copy of $\eta,$ namely, $\bar{\eta}\stackrel{d}{=}\eta$ and $\bar{\eta}$ is independent from $\eta.$
\\\textbf{Definition 2.5 (G-normal distribution).}  A $d$-dimensional
random vector $X=(X_1,\cdots,X_d)^T$ on a sublinear expectation
space $(\k)$ is called G-normal distribution
if$$aX+b\bar{X}\stackrel{d}{=}\sqrt{a^2+b^2}X,\ \ \forall a,b\geq 0,$$ where
$\bar{X}$ is an independent copy of $X.$

Given a pair of d-dimensional random vectors$(X,Y )$ in a sublinear
expectation space$(\k)$, Peng  defined a function
\begin{equation} G(P,A):=\e[\frac{1}{2}\langle
AX,X\rangle+\langle P,Y\rangle], \ \ (P,A)\in R^d\times S(d).
\end{equation}
It is easy to check that $G : R^d \times S(d)\rightarrow R$ is a
sublinear function monotonic in $A \in S(d)$ in the following sense:
for each $P,\bar{P}\in R^d$ and $A, \bar{A} \in S(d)$

\begin{equation} \label{eq:1}
\left\{ \begin{aligned}
         G(P+\bar{P},A+\bar{A}) &\leq G(P,A)+G(\bar{P},\bar{A}) ,\\
                  G(\lambda P,\lambda A)&=\lambda G(P,A),\forall \lambda\geq
                  0,\\
                  G(P,A)&\geq G(p,\bar{A}),if A\geq \bar{A}.
                          \end{aligned} \right.
                          \end{equation}
                          Clearly, G is also a continuous function.\\
\textbf{Proposition 2.1} (see Peng \cite{p6}). Let $G : R^d \times S(d) \rightarrow R$ be
a given sublinear and continuous function,  monotonic in $A\in S(d)$ in the sense of
(2.2). Then there exists a G-normal distributed $d$-dimensional random vector $X$ and a maximal distributed $d$-dimensional random vector $\eta$  on some sublinear
expectation space $(\k)$ satisfying (2.1) and
\begin{equation}(aX + b \bar{X},a^2 \eta + b^2 \bar{\eta} ) \stackrel{d}{=}(\sqrt{
a^2 + b^2}X, (a^2 + b^2 )\eta ), \forall a,b \geq 0,\end{equation} where
$( \bar{X},\bar{ \eta }) $ is an independent copy of $(X,\eta ).$
\\
\textbf{Definition 2.6 (G-distributed).} The pair $(X,\eta )$
satisfying (2.3) is called G-distributed.
\\
\textbf{Proposition 2.2} (see Peng \cite{p6}). For the pair $(X,\eta)$ satisfying (2.3)
and a function $\varphi\in C_{l,lip}(R^d\times R^d),$ then
$$u(t,x,y):=\e[\varphi(x+\sqrt{X},y+t\eta)],\ \  (t,x,y)\in [0,\infty)\times R^d\times R^d$$
is the unique viscosity solution of the PDE
\begin{equation}
\left\{ \begin{aligned}
         &\partial_tu-G(D_yu,D_x^2u)=0, \\
                  &u|_{t=0}=\varphi,
                          \end{aligned} \right.
                          \end{equation}
where $G:R^d\times
S(d)\rightarrow R$ is defined by (2.1) and
$D_yu=(\partial_{y_i}u)_{i=1}^d, \ D^2u=(\partial_{x_ix_j}^2u)_{i,j=1}^d.$
\section{Main results}
 \setcounter{equation}{0}
In this section, we give our main result and proof.
\\ \textbf{Theorem 3.1.}  Let $\{(X_i,Y_i)\}_{i=1}^{\infty}$ be a sequence of $R^d\times R^d$- valued random vectors on a sublinear expectation space $(\k)$ satisfying the following conditions:\\
(1) $(X_{i+1},Y_{i+1})$ is independent from
$\{(X_1,Y_1),\cdots,(X_i,Y_i)\}$ for each  $i=1,2,\cdots$; \\
(2) $\e[X_i]=\e[-X_i]=0,\
\e[X_{i}^2]=\bar{\sigma}_i^2,  \ -\e[-X_i^2]=\underline{\sigma}_i^2,$ where $0\leq \underline{\sigma}_i\leq \bar{\sigma}_i<\infty;$\\
(3) $\e[Y_i]=\bar{\mu}_i, \  -\e[-Y_i]=\underline{\mu}_i,\ \e[|X_{i }|^{2+\alpha}]+\e[|Y_{i }|^{2+\alpha}]<\infty
$ for every
$\alpha\in(0,1).$\\
(4) Let $\{w_i,i\geq 1\}$ be a
sequence of reals satisfying
$$\lim\limits_{n\rightarrow \infty}\frac{\sum\limits_{i=0}^{n-1}w_{i+1}^{2+\alpha}}{W_n^{1+\alpha/2}}=0,$$
where $W_n=\sum\limits_{i=1}^nw_i^2.$\\
(5) There are constants $\bar{\mu}\in R^d,\underline{\mu}\in R^d,\bar{\sigma}\in (R^+)^d,\underline{\sigma}\in (R^+)^d$ such that
$$\lim\limits_{n\rightarrow \infty}\sum_{i=0}^{n-1}\frac{w_{i+1}^2}{W_n}(|\bar{\sigma}_{i+1}^2-\bar{\sigma}^2|+|\underline{{\sigma}}^2-\underline{{\sigma}}^2_{i+1}|)=0,$$
and $$\lim\limits_{n\rightarrow \infty}\sum_{i=0}^{n-1}\frac{w_{i+1}^2}{W_n}(|\bar{\mu}_{i+1}-\bar{\mu}|+|\underline{{\mu}}-\underline{{\mu}}_{i+1}|)=0.$$
Then for every test function $\varphi\in
C_{b,lip}(R^d),$
 we have
 \begin{equation*}
\begin{split}
\lim\limits_{n\rightarrow\infty}\e[\varphi({{\bar{S}_n}})]=\e[\varphi(X+\eta)],
 \end{split}
 \end{equation*}
 where
$\{{{\bar{S}_n}}\}$ is defined by
$${{\bar{S}_n}}:=\sum\limits_{i=1}^n(\frac{w_iX_i}{\sqrt{W_n}}+\frac{w_i^2Y_i}{W_n}),$$  and the pair $(X,\eta)$
is  G-distributed.  The corresponding sublinear function
$G:R^d\times S(d)\rightarrow R$ is defined by
$$G(P,A):=\frac{1}{2}(\langle A^+\bar{\sigma},\bar{\sigma}\rangle-\langle A^-\underline{{\sigma}},\underline{{\sigma}}\rangle)+\langle P^+,\bar{\mu}\rangle-\langle P^-,\underline{{\mu}}\rangle,\ \ P\in R^d,A\in S(d),$$
where $A^+:=\max(A,\textbf{0}), A^-:=\max(-A,\textbf{0}).$
\\ \textbf{Proof.} For a small but fixed $h>0.$ Let
$V$ be the unique viscosity solution of
\begin{equation}
\partial_tV+G(D_xV,D_x^2V)=0,\ (t,x)\in [0,1+h)\times R^d,\ V|_{t=1+h}=\varphi.\end{equation}
From Proposition 2.2, we have
$$V(t,x)=\e[\varphi(x+\sqrt{1+h-t}X+(1+h-t)\eta)],$$
thus $$V(h,0)=\e[\varphi(X+\eta)],\ \ V(1+h,x)=\varphi(x).$$

 Because
\begin{equation}
\begin{split}
|\e[\varphi(\bar{{S_n}})]-\e[\varphi(X+\eta)]|&=|\e[V(1+h,\bar{{S_n}})-V(h,0)]| \\
 &=
 |\e[V(1+h,{{\bar{S}_n}})]-\e[V(1,{{\bar{S}_n}})]\\&+\e[V(1,{{\bar{S}_n}})]-V(0,0)+V(0,0)-V(h,0)]|\\
&\leq
|\e[V(1+h,{{\bar{S}_n}})]-\e[V(1,{{\bar{S}_n}})]|\\&+|V(0,0)-V(h,0)|+|\e[V(1,{{\bar{S}_n}})]-V(0,0)|,
 \end{split}
 \end{equation}
then we have
 \begin{equation}
\begin{split}
|\e[V(1,{{\bar{S}_n}})]-\e[V(1+h,{{\bar{S}_n}})]|&=|\e[\varphi(\bar{S}_n)]-\e[\varphi(\bar{S}_n+\sqrt{h}X+h
\eta)]|\\
&\leq C_1(\sqrt{h}\e[|X|]+h\e[|\eta|]),
 \end{split}
 \end{equation}
and
\begin{equation}
\begin{split}
|V(0,0)-V(h,0)|&=|\e[\varphi(\sqrt{1+h}X+(1+h) \eta)]-\e[\varphi(X+\eta)]|\\
&\leq C_2((\sqrt{1+h}-1)\e[|X|]+h\e[|\eta|]).
 \end{split}
 \end{equation}

 Next, let's compute the  bound of $|V(1,{{\bar{S}_n}})-V(0,0)|$, we set
${\bar{S}_0}=0, \delta_0=0, \delta_{i}=\frac{\sum\limits_{j=1}^{i}w_j^2}{W_n},1\leq i\leq n$. Then
\begin{equation}
\begin{split}
V(1,{{{\bar{S}}_n}})-V(0,0)&=\sum_{i=0}^{n-1}\{V(\delta_{i+1},{{\bar{S}_{i+1}}})-V(\delta_{i},{{\bar{S}_{i}}})\}\\
&=\sum_{i=0}^{n-1}\{[V(\delta_{i+1},{{\bar{S}_{i+1}}})-V(\delta_{i},{{\bar{S}_{i+1}}})]+[V(\delta_{i},{{\bar{S}_{i+1}}})-V(\delta_{i},{{\bar{S}_{i}}})]\}\\
&:=\sum_{i=0}^{n-1}\{A_i+B_i\},
 \end{split}
 \end{equation}
where
$$A_i=\frac{w_{i+1}^2}{W_n}\partial_tV(\delta_i,\bar{S}_i)+\frac{1}{2}\langle D^2V(\delta_i,\bar{S}_i)X_{i+1},X_{i+1}\rangle \frac{w_{i+1}^2}{W_n}+
\langle
DV(\delta_i,\bar{S}_i),X_{i+1}\frac{w_{i+1}}{\sqrt{W_n}}+\frac{w_{i+1}^2}{W_n}Y_{i+1}\rangle,$$
\begin{equation}
\begin{split}
B_i&=\frac{w_{i+1}^2}{W_n}\int_0^1[\partial_tV(\delta_{i}+\beta \frac{w_{i+1}^2}{W_n},\bar{S}_{i+1})-\partial_tV(\delta_i,\bar{S}_{i+1})]d\beta
\\&+
[\partial_tV(\delta_i,\bar{S}_{i+1})-\partial_tV(\delta_i,\bar{S}_{i})]\frac{w_{i+1}^2}{W_n}\\
&+\langle
D^2V(\delta_i,\bar{S}_{i})\frac{w_{i+1}}{\sqrt{W_n}}X_{i+1},\frac{w_{i+1}^2}{W_n}Y_{i+1}\rangle+\frac{1}{2}\langle
D^2V(\delta_i,\bar{S}_{i})\frac{w_{i+1}^2}{W_n}Y_{i+1},\frac{w_{i+1}^2}{W_n}Y_{i+1} \rangle\\
&+\int_0^1\int_0^1\langle C_{\beta
r}^i(X_{i+1}\frac{w_{i+1}}{\sqrt{W_n}}+Y_{i+1}\frac{w_{i+1}^2}{W_n}),X_{i+1}\frac{w_{i+1}}{\sqrt{W_n}}+Y_{i+1}\frac{w_{i+1}^2}{W_n})\rangle
rd\beta dr
 \end{split}
 \end{equation}
 with $$C_{\beta r}^i=D^2V(\delta_i,\bar{S}_i+r\beta(X_{i+1}\frac{w_{i+1}}{\sqrt{W_n}}+Y_{i+1}\frac{w_{i+1}^2}{W_n}))-D^2V(\delta_i,\bar{S}_i).$$

 We note that
 $$\e[\langle
DV(\delta_{i},\bar{S}_i),X_{i+1}\frac{w_{i+1}}{\sqrt{W_n}}\rangle]=\e[\e[\langle\theta,X_{i+1}\frac{w_{i+1}}{\sqrt{W_n} }\rangle]_{\theta=DV(\delta_{i},\bar{S}_i)}]=0,$$ then
 \begin{equation}
 \begin{split}
 \e[A_i]&=\e[\frac{w_{i+1}^2}{W_n}\partial_tV(\delta_i,\bar{S}_i)+\frac{1}{2}\langle D^2V(\delta_i,\bar{S}_i)X_{i+1},X_{i+1}\rangle \frac{w_{i+1}^2}{W_n}+
\langle
DV(\delta_i,\bar{S}_i),\frac{w_{i+1}^2}{W_n}Y_{i+1}\rangle]\\
&\leq\frac{1}{2}\frac{w_{i+1}^2}{W_n}\e[(D^2V(\delta_i,\bar{S}_i))^+(\bar{\sigma}_{i+1}^2-\bar{\sigma}^2)+(D^2V(\delta_i,\bar{S}_i))^-(\underline{{\sigma}}^2-\underline{{\sigma}}^2_{i+1})]\\
&+\frac{w_{i+1}^2}{W_n}\e[(DV(\delta_i,\bar{S}_i))^+(\bar{\mu}_{i+1}-\bar{\mu})+(DV(\delta_i,\bar{S}_i))^-(\underline{{\mu}}-\underline{{\mu}}_{i+1})]\\
&\leq\frac{1}{2}\frac{w_{i+1}^2}{W_n}(|\bar{\sigma}_{i+1}^2-\bar{\sigma}^2|+|\underline{{\sigma}}^2-\underline{{\sigma}}^2_{i+1}|)\e[|D^2V(\delta_i,\bar{S}_i)|]\\
&+\frac{w_{i+1}^2}{W_n}(|\bar{\mu}_{i+1}-\bar{\mu}|+|\underline{{\mu}}-\underline{{\mu}}_{i+1}|)\e[|DV(\delta_i,\bar{S}_i)|].
  \end{split}
 \end{equation}

It is easy to check that there exist constants $C_3, C_4 > 0,$ such that $\e[|D^2V(\delta_i,\bar{S}_i)|]\leq C_3$ and $\e[|DV(\delta_i,\bar{S}_i)|]\leq C_4.$  Then we obtain
 \begin{equation}
 \begin{split}
 \e[\sum_{i=0}^{n-1}A_i]&\leq \sum_{i=0}^{n-1}\e[A_i]\\
 &\leq \frac{C_3}{2}\sum_{i=0}^{n-1}\frac{w_{i+1}^2}{W_n}(|\bar{\sigma}_{i+1}^2-\bar{\sigma}^2|+|\underline{{\sigma}}^2-\underline{{\sigma}}^2_{i+1}|)\\
 &+C_4\sum_{i=0}^{n-1}\frac{w_{i+1}^2}{W_n}(|\bar{\mu}_{i+1}-\bar{\mu}|+|\underline{{\mu}}-\underline{{\mu}}_{i+1}|).
 \end{split}
 \end{equation}

By (5) in Theorem 3.1, we have \begin{equation}
 \begin{split}
 \lim\limits_{n\rightarrow \infty}\e[\sum_{i=0}^{n-1}A_i]\leq 0.
  \end{split}
 \end{equation}

 Similarly, we also have
 \begin{equation}
 \begin{split}
 \lim\limits_{n\rightarrow \infty}\e[\sum_{i=0}^{n-1}A_i]\geq 0.
  \end{split}
 \end{equation}

 Thus
 \begin{equation}
 \begin{split}
 \lim\limits_{n\rightarrow \infty}\e[\sum_{i=0}^{n-1}A_i]= 0.
  \end{split}
 \end{equation}

 For $B_i,$ since both $\partial_tV$ and $D^2V$ are uniformly
$\frac{\alpha}{2}$-h\"{o}lder continuous in $t$ and
$\alpha$-h\"{o}lder continuous in $x$ on $[0,1]\times R^d,$ then
\begin{equation}
\begin{split}
\e[|B_i|]&\leq
C_5(\frac{w_{i+1}^2}{W_n})^{1+\frac{\alpha}{2}}+C_6\e[|X_{i+1}\frac{w_{i+1}}{\sqrt{W_n}}+\frac{w_{i+1}^2}{W_n}Y_{i+1}|^{\alpha}]\frac{w_{i+1}^2}{W_n}
\\&+C_7\e[|X_{i+1}\frac{w_{i+1}}{\sqrt{W_n}}+\frac{w_{i+1}^2}{W_n}Y_{i+1}|^{2+\alpha}]\\
&\leq
C_8\frac{w_{i+1}^{2+\alpha}}{W_n^{1+\alpha/2}}(1+\e[|X_{i+1}|^{\alpha}]+\e[|X_{i+1}|^{2+\alpha}+\e[|Y_{i+1}|^\alpha]+\e[|Y_{i+1}|^{2+\alpha}])\\
&\leq C_9\frac{w_{i+1}^{2+\alpha}}{W_n^{1+\alpha/2}}.
 \end{split}
 \end{equation}
 Thus
 \begin{equation}
\begin{split}
-C_9\frac{\sum\limits_{i=0}^{n-1}w_{i+1}^{2+\alpha}}{W_n^{1+\alpha/2}}&\leq \e[V(1,\bar{S}_n)]-V(0,0)\leq C_9\frac{\sum\limits_{i=0}^{n-1}w_{i+1}^{2+\alpha}}{W_n^{1+\alpha/2}}.
 \end{split}
 \end{equation}
By the condition (4) in Theorem 3.1,  as $n\rightarrow \infty$,  we have  \begin{equation}
\begin{split}
\lim\limits_{n\rightarrow \infty}\e[V(1,\bar{S}_n)]=V(0,0).
 \end{split}
 \end{equation}

 It follows from (3.2)-(3.4) and (3.14) that
 \begin{equation}
\begin{split}
\limsup\limits_{n\rightarrow\infty}|\e[\varphi(\bar{S}_n)]-\e[\varphi(X+\eta)]|\leq C_{10}(h+\sqrt{h}).
 \end{split}
 \end{equation}
Since $h$ can be arbitrarily small, we have,
$$\lim\limits_{n\rightarrow\infty}\e[\varphi(\bar{S}_n)]=\e[\varphi(X+\eta)].$$
The proof is complete.

From the Theorem 3.1,  we can obtain  the following corollaries:
\\\textbf{Corollary 3.1.}  Let $\{X_i\}_{i=1}^{\infty}$ be a sequence of $R^d$- valued random vectors on a sublinear expectation space $(\k)$ satisfying the following conditions:\\
(1) $X_{i+1}$ is independent from
$(X_1,\cdots, X_i)$ for each  $i=1,2,\cdots$; \\
(2) $\e[X_i]=\e[-X_i]=0,\
\e[X_{i}^2]=\bar{\sigma}_i^2,  \ -\e[-X_i^2]=\underline{\sigma}_i^2,$ where $0\leq \underline{\sigma}_i\leq \bar{\sigma}_i<\infty;$\\
(3) $ \e[|X_{i }|^{2+\alpha}]<\infty
$ for every
$\alpha\in(0,1).$\\
(4) Let $\{w_i,i\geq 1\}$ be a
sequence of reals satisfying
$$\lim\limits_{n\rightarrow \infty}\frac{\sum\limits_{i=0}^{n-1}w_{i+1}^{2+\alpha}}{W_n^{1+\alpha/2}}=0,$$
where $W_n=\sum\limits_{i=1}^nw_i^2.$\\
(5) There are constants $\bar{\sigma}\in (R^+)^d,\underline{\sigma}\in (R^+)^d$ such that
$$\lim\limits_{n\rightarrow \infty}\sum_{i=0}^{n-1}\frac{w_{i+1}^2}{W_n}(|\bar{\sigma}_{i+1}^2-\bar{\sigma}^2|+|\underline{{\sigma}}^2-\underline{{\sigma}}^2_{i+1}|)=0.$$
Then for every test function $\varphi\in
C_{b,lip}(R^d),$
 we have
 \begin{equation*}
\begin{split}
\lim\limits_{n\rightarrow\infty}\e[\varphi({{\bar{S}_n}})]=\e[\varphi(X)],
 \end{split}
 \end{equation*}
 where
$\{{{\bar{S}_n}}\}$ is defined by
$${{\bar{S}_n}}:=\sum\limits_{i=1}^n\frac{w_iX_i}{\sqrt{W_n}},$$  and $X$
is a G-normal distributed random vector  and the corresponding sublinear function
$G: S(d)\rightarrow R$ is defined by
$$G(A):=\frac{1}{2}(\langle A^+\bar{\sigma},\bar{\sigma}\rangle-\langle A^-\underline{{\sigma}},\underline{{\sigma}}\rangle),\ \ A\in S(d).$$
\textbf{Remark 3.1.} Let $w_i\equiv 1,i=1,\cdots,$ and $d=1$ in Corollary 3.1, then it is the main result of  Li and Shi (see \cite{l}).
\\ \textbf{Corollary 3.2.}   Let $\{(X_i,Y_i)\}_{i=1}^{\infty}$ be a sequence of $R^d\times R^d$- valued random vectors on a sublinear expectation space $(\k)$ satisfying the following conditions:\\
(1) $(X_{i+1}, Y_{i+1})\stackrel{d}{=}(X_{i},Y_{i})$ for each  $i=1,2,\cdots$;\\
(2) $(X_{i+1},Y_{i+1})$ is independent from
$\{(X_1,Y_1),\cdots,(X_i,Y_i)\}$ for each  $i=1,2,\cdots$; \\
(3) $\e[X_1]=\e[-X_1]=0,\
\e[|X_{1 }|^{2+\alpha}]+\e[|Y_{1 }|^{2+\alpha}]<\infty
$ for every
$\alpha\in(0,1).$\\
(4) Let $\{w_i,i\geq 1\}$ be a
sequence of reals satisfying
$$\lim\limits_{n\rightarrow \infty}\frac{\sum\limits_{i=0}^{n-1}w_{i+1}^{2+\alpha}}{W_n^{1+\alpha/2}}=0,$$
where $W_n=\sum\limits_{i=1}^nw_i^2.$\\
Then for every test function $\varphi\in
C_{b,lip}(R^d),$
 we have
 \begin{equation*}
\begin{split}
\lim\limits_{n\rightarrow\infty}\e[\varphi({{\bar{S}_n}})]=\e[\varphi(X+\eta)],
 \end{split}
 \end{equation*}
 where
$\{{{\bar{S}_n}}\}$ is defined by
$${{\bar{S}_n}}:=\sum\limits_{i=1}^n(\frac{w_{i}X_i}{\sqrt{W_n}}+\frac{w_i^2Y_i}{W_n}),$$  and the pair $(X,\eta)$
is G-distributed  and the corresponding sublinear function
$G:R^d\times S(d)\rightarrow R$ is defined by
$$G(P,A):=\e[\frac{1}{2}\langle AX_1,X_1\rangle+\langle P,Y_1\rangle],\ \ P\in R^d,A\in S(d).$$
\textbf{Remark 3.2.}  Let $w_i\equiv 1,i=1,\cdots,$  in Corollary 3.2, we can get the result of  Peng (see \cite{p3}).
\textbf{Remark 3.3.}  Let $w_i\equiv 1,i=1,\cdots,$  in Corollary 3.2, we can get the result of  Peng (see \cite{p3}).
\\ \textbf{Corollary 3.3.}   Let $\{Y_i\}_{i=1}^{\infty}$ be a sequence of $R^d$- valued random vectors on a sublinear expectation space $(\k)$ satisfying the following conditions:\\
(1) $Y_{i+1}$ is independent from
$\{Y_1,\cdots,Y_i\}$ for each  $i=1,2,\cdots$; \\
(2) $\e[Y_i]=\bar{\mu}_i, \  -\e[-Y_i]=\underline{\mu}_i,\ \e[|Y_{i }|^{2+\alpha}]<\infty
$ for every
$\alpha\in(0,1).$\\
(3) We assume that there exist $\bar{\mu}$ and  ${\underline{\mu}}$ such that
$$\lim\limits_{n\rightarrow \infty}\frac{1}{n}\sum_{i=0}^{n-1}|\bar{\mu}-\bar{\mu}_{i+1}|=0,$$
and
$$\lim\limits_{n\rightarrow \infty}\frac{1}{n}\sum_{i=0}^{n-1}|{\underline{\mu}}_{i+1}-{\underline{\mu}}|=0.$$
Then for every test function $\varphi\in
C_{b,lip}(R^d),$
 we have
 \begin{equation*}
\begin{split}
\lim\limits_{n\rightarrow\infty}\e[\varphi({{\frac{\sum\limits_{i=1}^nY_i}{n}}})]=\e[\varphi(\eta)],
 \end{split}
 \end{equation*}
 where
 $\eta$
is  maximal distribution.  The corresponding sublinear function
$G:R^d\rightarrow R$ is defined by
$$G(P):=\langle P^+,\bar{\mu}\rangle-\langle P^-,\underline{{\mu}}\rangle,\ \ P\in R^d.$$
\textbf{Remark 3.4.} If $\bar{\mu}_i=\bar{\mu}$ and $\underline{{\mu}}_i=\underline{{\mu}}$ for all $1\leq i\leq n$, then it is the law of large numbers under sublinear space  (see Peng \cite{p3}).
\section*{Acknowledgements}
The author  thanks  Prof.  Zengjing Chen for helpful discussion and valuable suggestions.

The author thanks  the partial support from The National Basic Research Program of China (973 Program)(No. 2007CB814901).

\end{document}